\documentclass[11pt]{amsart}
\usepackage{amsfonts}
\usepackage{amsmath}
\usepackage{amssymb}
\usepackage{amsthm}
\usepackage{graphicx}
\usepackage{enumerate}
\usepackage{multicol}
\usepackage{mathrsfs}
\usepackage{xypic}
\usepackage{enumerate}

\DeclareFontFamily{U}{mathx}{\hyphenchar\font45}
\DeclareFontShape{U}{mathx}{m}{n}{
      <5> <6> <7> <8> <9> <10>
      <10.95> <12> <14.4> <17.28> <20.74> <24.88>
      mathx10
      }{}
\DeclareSymbolFont{mathx}{U}{mathx}{m}{n}
\DeclareFontSubstitution{U}{mathx}{m}{n}
\DeclareMathAccent{\widecheck}{0}{mathx}{"71}
\DeclareMathAccent{\wideparen}{0}{mathx}{"75}

\newcounter{dummy} \numberwithin{dummy}{section}
\newtheorem{theorem}[dummy]{Theorem}

\newtheorem{lemma}[dummy]{Lemma}

\newtheorem{proposition}[dummy]{Proposition}
\theoremstyle{remark}
\newtheorem{remark}[dummy]{Remark}
\newtheorem{example}[dummy]{Example}

\numberwithin{equation}{section}

\newcommand{\calE}{\mathcal E}
\newcommand{\calF}{\mathcal F}
\newcommand{\calH}{\mathcal H}
\newcommand{\calL}{\mathcal L}
\newcommand{\calR}{\mathcal R}
\newcommand{\calV}{\mathcal V}

\DeclareMathOperator{\pr}{pr}
\DeclareMathOperator{\Sym}{Sym}
\DeclareMathOperator{\Ann}{Ann}

\DeclareMathOperator{\tensorh}{\mathbf h}
\DeclareMathOperator{\tensorg}{\mathbf g}

\DeclareMathOperator{\tensorv}{\mathbf v}
\DeclareMathOperator{\II}{II}
\DeclareMathOperator{\vl}{vl}

\DeclareMathOperator{\rnabla}{\mathring{\nabla}}
\DeclareMathOperator{\tensors}{\mathbf{s}}
\DeclareMathOperator{\shh}{\sharp^{{\bf h}^*}}
\DeclareMathOperator{\shv}{\sharp^{{\bf v}^*}}
\DeclareMathOperator{\shs}{\sharp^{{\bf s}^*}}

\title{Riemannian and Sub-Riemannian geodesic flows}
\author{Mauricio Godoy Molina \\ Erlend Grong}

\address{Department of Mathematics, University of Bergen, Norway \&~Departamento de Matem\'aticas. Universidad de la Frontera, Chile.}
\email{mauricio.godoy@math.uni.lu}

\address{Mathematics Research Unit, University of Luxembourg, Luxembourg.}
\email{erlend.grong@uni.lu}



\begin{document}

\begin{abstract}
In the present paper we show that the geodesic flows of a sub-Riemannian metric and that of a Riemannian extension commute if and only if the extended metric is parallel with respect to a certain connection. This helps us to describe the geodesic flow of sub-Riemannian metrics on totally geodesic Riemannian submersions. As a consequence we can characterize sub-Riemannian geodesics as the horizontal lifts of projections of Riemannian geodesics.
\end{abstract}

\maketitle

\section{Introduction}
Since the introduction of sub-Riemannian geometry in 1986, see \cite{Str86}, one of the main topics of research has been finding important geometric invariants. A sub-Riemannian manifold is a connected manifold $M$ with a smoothly varying inner product $\tensorh$ defined only on a subbundle $\calH$ of the tangent bundle. Such spaces have a metric structure by considering the distance between two points as the infimum of the length of all curves tangent to $\calH$ connecting them. Understanding what the proper generalization of curvature should be for such spaces, has been a topic of great interest in recent years.

There are currently two main ways of attacking this problem. One approach considers symplectic invariants of the (normal) geodesic flow on sub-Riemannian spaces. For results in this direction, see e.g.~\cite{ABR13,BaLu14,LiZe11,ZeLi07}. The second approach tries to understand curvature in terms of properties of the heat flow corresponding to a second order differential operator, known as sub-Laplacian, which is the sub-Riemannian analogue of the Laplace-Beltrami operator. Unlike the first mentioned approach, the second one requires a extension of $\tensorh$ to a Riemannian metric satisfying certain properties. Such an extension will not be unique in general, see \cite[Section~4.5]{GrTh14a}. However, once an appropriate choice has been made, one has at hand powerful results such as a parabolic Harnack-inequality and a Bonnet-Myers theorem, see e.g. \cite{BaGa12,BBG14,BKW14,GrTh14b}.

In the present paper we give an initial step to bring together both approaches, by studying the Riemannian and sub-Riemannian geodesic flows. Our main result states that these flows commute if and only if a Riemannian metric taming $\tensorh$ is parallel with respect to an appropriate connection. This requirement also appears as a hypothesis in the second approach to curvature described above, see Remark~\ref{rem:2ndapproach}. Moreover, we prove that requiring that the projection of the Riemannian and sub-Riemannian geodesic flows to the base space of a submersion coincide is equivalent to requiring the fibers of the submersion to be totally geodesic. This generalizes a result found in \cite{Mon84}, where it was shown that the trajectories of particles with a given gauge in a Yang-Mills field can be considered as projections of both sub-Riemannian and Riemannian geodesics on a principal bundle, given that the gauge group has a bi-invariant metric. See Example~\ref{ex:Montgomery} for more details.

The structure of the paper is as follows. In Section~\ref{sec:SRR} we introduce the main concepts that we will use and state our results. We postpone the proofs to the next section, for the sake of clarity. Section~\ref{sec:SRR} concludes with two relevant examples. The technical tools and the proofs of our results are presented in Section~\ref{sec:proofs}.

\subsection{Acknowledgments} 
We thank Dr. Petri Kokkonen for helpful discussions and comments.

\subsection{Notation and conventions}
All manifolds are smooth and connected.
For any vector bundle $\calE \to M$ over a manifold $M$, we will use $\Gamma(\calE)$ for the space of all smooth sections of $\calE$. For any vector field $X \in \Gamma(TM)$, we will write $\calL_X$ for the Lie derivative with respect to $X$ and $e^{tX}$ for its local flow on $M$. If $\calE$ is a subbundle of the tangent bundle $TM$, then $\Ann(\calE)$ denotes the subbundle of $T^*M$ of all covectors that vanish on $\calE$.
\section{Geodesic flows: Statement of the results} \label{sec:SRR}
\subsection{Sub-Riemannian manifolds}
\emph{A sub-Riemannian manifold} is a triple $(M, \calH, \tensorh)$, where $M$ is a (connected) manifold, $\calH$ is a subbundle of the tangent bundle $TM$ and $\tensorh$ is a metric tensor defined only on $\calH$. Equivalently, it can be considered as a pair $(M, \tensorh^*)$, where $M$ is a manifold and $\tensorh^*$ is a bilinear positive semidefinite tensor of the cotangent bundle that vanishes on a subbundle of $T^*M$. We will call $\tensorh^*$ a sub-Riemannian \emph{cometric}. The relation between $(\calH,\tensorh)$ and $\tensorh^*$ can be described as follows. Let $\calH$ be the image of the map $\shh$ given by
\begin{equation} \label{cor1} \shh: T^*M \to TM, \qquad p \mapsto \tensorh^*(p, \, \centerdot \,).\end{equation}
and endow $\calH$ with a metric tensor $\tensorh$ determined by equation
\begin{equation} \label{cor2} \tensorh(\shh p_1, \shh p_2) := \tensorh^*(p_1, p_2), \qquad p_j \in T_xM, x \in M, j=1,2.\end{equation}
Conversely, given the pair $(\calH, \tensorh)$, the cometric $\tensorh^*$ is uniquely determined by \eqref{cor1} and \eqref{cor2}. The kernel of $\shh$ will be the subbundle $\Ann(\calH)$ of $T^*M$.

An absolutely continuous curve in $(M, \calH, \tensorh)$ is called \emph{horizontal} if $\dot \gamma(t) \in \calH_{\gamma(t)}$ for almost every $t$. The distance in a sub-Riemannian manifold is given by \emph{the Carnot-Carath\'eodory metric} $\mathsf{d}^{\tensorh}$ defined so that $\mathsf{d}^{\tensorh}(x,y)$ is the infimum of all integrals $\int_0^1 \tensorh(\dot \gamma(t), \dot \gamma(t) )^{1/2} \, dt$ taken over all horizontal curves $\gamma$ satisfying $\gamma(0) = x$ and $\gamma(1) = y.$ This distance can only be finite if any two points can be connected by a horizontal curve. A sufficient condition for this to hold, is that $\calH$ is \emph{bracket-generating}, i.e. that the sections of $\calH$ and their iterated brackets span $TM$.

Minimizers of the distance $\mathsf{d}^{\tensorh}$ are either \emph{normal geodesics} or \emph{abnormal curves}. Normal geodesics are projections of integral curves of the Hamiltonian vector field $\vec{H}^{\tensorh}$ of the Hamiltonian $H^{\tensorh}(p) := \frac{1}{2} \tensorh^*(p,p).$ Such curves are always locally length minimizers and smooth. For the definition of abnormal curves and more details on sub-Riemannian manifolds in general, we refer to~\cite{Mon02}.

For future computations, the following description of normal geodesics with respect to an arbitrary affine connection $\nabla$ will be convenient.
\begin{proposition} \label{prop:NormGeo}
Let $\nabla$ be any affine connection on $M$ with torsion $T^{\nabla}$. Then a curve $\lambda(t)$ in $T^*M$ with projection $\gamma$ is an integral curve of $\vec{H}^{\tensorh}$ if and only if
$$\dot \gamma(t) = \shh \lambda(t), \qquad \nabla_{\dot \gamma} \lambda(t) = -\lambda(t) T^{\nabla}(\dot \gamma, \, \centerdot \, ) + (\nabla_{\, \centerdot \,} \tensorh^*)(\lambda(t), \lambda(t)).$$
\end{proposition}
Let $\Pi^M:T^*M \to M$ be the canonical projection of the cotangent bundle.
For every $x \in M$, following \cite{Str86}, we define the sub-Riemannian exponential $\exp^{sr}(x, \, \centerdot \,): U_x \subseteq T^*M \to M$ by the formula
\begin{equation} \label{expsr}\exp^{sr}(x,p ) : = \lambda(1), \qquad \lambda(t) := (\Pi^M \circ e^{t\vec{H}^{\tensorh}})(p),\end{equation}
where $U_x$ is the collection of all $p \in T^*_xM$ such that \eqref{expsr} is well defined.

\subsection{Taming sub-Riemannian metrics}
Let $(M, \calH, \tensorh)$ be a a sub-Riemannian manifold and let $\tensorg$ denote a Riemannian metric on $M$ such that $\tensorg|_{\calH} = \tensorh$. Such a Riemannian metric $\tensorg$ is said to \emph{tame} $\tensorh$.  Let $\nabla^{\tensorg}$ denote the Levi-Civita connection associated to $\tensorg$ and define $\calV$ as the orthogonal complement of $\calH$ with respect to $\tensorg$. We will use $\pr_{\calH}$ and $\pr_{\calV}$ for the respective orthogonal projections to $\calH$ and $\calV$. We introduce a connection $\rnabla$, which will play a central role in our results, as follows
\begin{align} \label{rnabla}
\rnabla_XY & := \pr_{\calH} \nabla_{\pr_{\calH} X}^{\tensorg} \pr_{\calH} Y + \pr_{\calV} \nabla_{\pr_{\calV} X}^{\tensorg} \pr_{\calV} Y \\ \nonumber
& \qquad + \pr_{\calH} [\pr_{\calV} X, \pr_{\calH} Y] + \pr_{\calV} [\pr_{\calH} X, \pr_{\calV} Y].
\end{align}
Define a metric tensor $\tensorv$ on $\calV$ by $\tensorv := \tensorg|_{\calV}$. This corresponds to a (degenerate) cometric $\tensorv^*$ on the cotangent bundle through the relations~\eqref{cor1} and~\eqref{cor2}. This cometric defines a Hamiltonian function $H^{\tensorv}(p) := \frac{1}{2} \tensorv^*(p,p) $ for any $ p \in T^*M$. Write $H^{\tensorg} = H^{\tensorh} + H^{\tensorv}$, which is the Hamiltonian of the Riemannian metric $\tensorg$. The following result relates the connection $\rnabla$ with the Hamiltonian functions defined by the cometrics $\tensorh^*$ and $\tensorv^*$.

\begin{lemma} \label{lemma:commute}
Let $\{\, \centerdot\, ,\, \centerdot \, \}$ denote the Poisson bracket with respect to the canonical symplectic form on $T^*M$. Then $\{ H^{\tensorh}, H^{\tensorv}\} = \{ H^{\tensorh}, H^{\tensorg}\} = 0$ if and only if $\rnabla \tensorg = 0$.
\end{lemma}
Let $\exp^{r}(x,\, \centerdot \,):V_x \subseteq T_xM \to M$ denote the Riemannian exponential map with respect to $\tensorg$ from the point $x$. Let $\sharp$ be the identification of $T^*M$ with $TM$ using $\tensorg$. It is clear that $\exp^r(x, t\sharp p) = \big(\Pi^M \circ e^{t\vec{H}^{\tensorg}}\big)(p)$, since projections of the solutions to the Hamiltonian system with Hamiltonian $H^{\tensorg}(p) = \frac{1}{2}\tensorg(\sharp p, \sharp p)$ are exactly the Riemannian geodesics. It follows from Lemma~\ref{lemma:commute} that if $\rnabla \tensorg =0$, then for $p \in T^*M$, we have
$$\exp^r(x, t\sharp p) = \big(\Pi^M \circ e^{t\vec{H}^{\tensorh}} \circ e^{t\vec{H}^{\tensorv}}\big)(p) = \big(\Pi^M \circ e^{t\vec{H}^{\tensorv}} \circ e^{t\vec{H}^{\tensorh}}\big)(p),$$
for any value of $t$ such that the above terms are well defined. Furthermore, $ e^{s\vec{H}^{\tensorh}} \circ e^{t\vec{H}^{\tensorg}}(p) =  e^{t\vec{H}^{\tensorg}} \circ e^{s\vec{H}^{\tensorh}}(p)$ for any $p\in T^*M$, and $t,s\in{\mathbb R}$ such that both sides are well defined.

\begin{remark}\label{rem:2ndapproach}
The connection $\rnabla$ was first introduced in \cite{BKW14,GrTh14a,GrTh14b} as a tool to obtain generalized curvature-dimension inequalities for sub-Riemannian manifolds. Such inequalities connect a Riemannian metric $\tensorg$ taming $\tensorh$ with the second order operator $\Delta^{\tensorh}$ given by
$$\Delta^{\tensorh} f = \mathrm{div} \, \shh df, \quad f \in C^\infty(M),$$
where the divergence is with respect to the volume form of $M$ defined by $\tensorg$. It turns out that the possibility of choosing a Riemannian extension such that $\rnabla \tensorg = 0$ is essential for obtaining results such as a parabolic Harnack inequality for the heat flow of $\Delta^{\tensorh}$. However, note that in this setting we also need the requirement that the trace of the map
\begin{equation} \label{requirement} 
X \mapsto \pr_{\calH} [ \pr_{\calV} Y, \pr_{\calV} [ \pr_{\calH} Y, \pr_{\calH}X]] ,
\end{equation} 
vanishes for any vector field $Y$. Although the map~\eqref{requirement} is written with vector fields $X$ and $Y$, it is in fact tensorial in both arguments. This latter requirement is needed to ensure that $\Delta^{\tensorh}$ commutes with the Laplace-Beltrami operator of $\tensorg$. See \cite[Appendix~A]{GrTh14b} for details. Note that \eqref{requirement} vanishes identically if $\calV$ is integrable, which is the case considered in Section~\ref{sec:TotallyR}.
\end{remark}

\begin{remark}
All the results in this paper are still valid if we only require that $\tensorv$ is a nondegenerate metric tensor on $\calV$. The same is also true if we consider pseudo sub-Riemannian metrics on $\calH$, assuming only that $\tensorh$ is nondegenerate.
\end{remark}

\subsection{Totally geodesic Riemannian foliations} \label{sec:TotallyR}
Let $(M, \calH, \tensorh)$ be a sub-Riemannian manifold and let $\tensorg$ be a Riemannian metric taming $\tensorh$. Assume that the orthogonal complement $\calV$ of $\calH$ with respect to $\tensorg$ is integrable, i.e. we assume that for any pair of vector fields $Z, W \in \Gamma(\calV)$ we have $[Z,W]|_x \in \calV_x, x \in M$. By the Frobenius theorem, we know that there exists a foliation $\calF$ of $M$ with leaves tangent to $\calV$.

Define $\rnabla$ as in \eqref{rnabla}. It is simple to verify that for any $v,w \in T_xM, x \in M,$
\begin{equation} \label{CovDerg} (\rnabla_{v} \tensorg)(w,w)= (\rnabla_{\pr_{\calH} v} \tensorg)(\pr_{\calV} w, \pr_{\calV} w) + (\rnabla_{\pr_{\calV} v} \tensorg)(\pr_{\calH} w, \pr_{\calH} w)
\end{equation}
It follows that $\rnabla \tensorg = 0$ if and only if for any $X \in \Gamma(\calH)$ and $Z \in \Gamma(\calV)$, we have
\begin{eqnarray} 
 \label{TGF} (\rnabla_X \tensorg)(Z,Z) &=& (\calL_X \tensorg)(Z,Z) = 0, \quad \text{and} \\
\label{RiemannF} (\rnabla_Z \tensorg)(X,X) &=& (\calL_Z \tensorg)(X,X) = 0. 
 \end{eqnarray}
To explain the geometric meaning of \eqref{TGF}, let $\II$ be \emph{the second fundamental form} of the leaves of the foliation $\calF$. This is a symmetric, bilinear vector valued tensor defined by $\II(Z,W) := \pr_{\calH} \nabla^{\tensorg}_{\pr_{\calV} Z} \pr_{\calV} W.$ If $\II=0$, then $\calF$ is called \emph{totally geodesic}. This means that the leaves of $\calF$ are totally geodesic submanifolds of $M$, i.e. if $\calF_x$ is the leaf of $\calF$ containing $x \in M$ and $v \in T_x \calF_x = \calV_x \subseteq T_xM$, then the curve $\gamma(t) = \exp(x,tv)$ is contained in $\calF_x$. Using the definition of the Levi-Civita connection, it follows that $(\calL_X \tensorg)(Z,Z) = - 2 \tensorg(X, \II(Z,Z))$, hence, \eqref{TGF} is equivalent to $\calF$ being a totally geodesic foliation.

To understand the meaning of \eqref{RiemannF}, let us first start with the definition of a Riemannian submersion. Let $\pi: (M, \tensorg) \to (B, \widecheck{\tensorg})$ be a submersion between two Riemannian manifolds with $\calV = \ker \pi_*$ and $\calH = \calV^\perp$. The submersion $\pi$ is called \emph{Riemannian} if
$$\tensorg(v, w) = \widecheck{\tensorg}(\pi_* v, \pi_* w), \qquad \text{for any } v,w \in \calH.$$
Let $\calF$ be a  foliation on a Riemannian manifold $(M, \tensorg)$ satisfying \eqref{RiemannF}. Then for every $x \in M$, there exists a neighborhood $U$ of $x$ such that the quotient $B= U/\calF|_U$ is a well defined manifold that can be given a metric $\widecheck{\tensorg}$ such that the quotient map $\pi:(U, \tensorg|_U) \to (B, \widecheck{\tensorg})$ is a Riemannian submersion. Hence, we call such foliations \emph{Riemannian}, since thay can locally be obtained from a Riemannian submersion, see e.g.~\cite{Rei59} for more details.

Using Lemma~\ref{lemma:commute}, we have the following result
\begin{theorem} \label{th:TGRF}
Let $(M, \tensorg)$ be a Riemannian manifold and let $\calF$ be a totally geodesic Riemannian foliation of $M$ corresponding to an integrable subbundle~$\calV$. Define a sub-Riemannian manifold $(M, \calH, \tensorh)$ where $\calH$ is the orthogonal complement of $\calV$ and $\tensorh = \tensorg|_{\calH}$. Then, for any $x \in M$ and $p \in T^*M$, if $P_t$ denotes parallel transport of vectors along $\exp^r(x, t \sharp p)$ with respect to the Levi-Civita connection, we have
$$\exp^{sr}(x,tp) = \exp^r\left(\exp^r(x,t\sharp p),- t \pr_{\calV} P_t \sharp p\right),$$
for any $t$ such that both sides are well defined.
\end{theorem}

\begin{remark}
It follows from the proof of this theorem, found in Section~\ref{sec:ProofTGRF}, that we could have also written
$$\exp^{sr}(x,tp) = \exp^r\left(\exp^r(x, - t\pr_{\calV} \sharp p),t \widetilde P_t \sharp p\right),$$
where $\widetilde P_t$ denotes the parallel transport along $\exp^r(x, - \pr_{\calV} t\sharp p)$ with respect to the Levi-Civita connection.
\end{remark}

\subsection{Submersions and sub-Riemannian geometry} \label{sec:submersion}
Let $(B, \widecheck{\tensorg})$ be a Riemannian manifold and let $\pi: M \to B$ be a submersion into $B$ with vertical bundle $\calV := \ker \pi_*$. Since the vector fields with values in $\calV$ are exactly the vector fields on $M$ that are $\pi$-related to the zero section of $TB$, $\calV$ is an integrable subbundle. The leaves of the corresponding foliation are given by submanifolds $M_b := \pi^{-1}(b), b \in B$. Moreover, $M_b$ is an embedded submanifold since any $b \in B$ is a regular value of $\pi$.

\emph{An Ehresmann connection} $\calH$ on $\pi$ is a subbundle satisfying $TM = \calH \oplus \calV$. Each Ehresmann connection $\calH$ on $\pi$ gives us a sub-Riemannian structure $(\calH, \tensorh)$ by lifting $\widecheck{\tensorg}$, i.e.
$$\tensorh(v,w) = \widecheck{\tensorg}(\pi_* v, \pi_* w), \qquad v, w \in \calH_x, x \in M.$$
Let $\tensorg$ be any Riemannian metric on $M$ taming $\tensorh$. We are interested in a ``good'' way of choosing the Riemannian metric $\tensorg$, in the sense that we want to consider when $\tensorg$ satisfies the property 
\begin{equation} \label{CondA}
\tag{${\sf A}$} \begin{array}{l} \text{For any $x \in M$, there is a neighborhood $U_x$ of $0 \in T^*_xM$ such}\\
\text{that for any $p \in U_x $, the curves $\gamma(t) = \exp^{sr}(x, tp)$ and} \\
\text{ $\eta(t) = \exp^r(x,t \sharp p)$, $0 \leq t \leq 1$, have the same projection in $M$. } \end{array}
\end{equation}

Related to this condition, we have the following result for submersions.
\begin{theorem} \label{th:A}
The condition \eqref{CondA} holds if and only if
\begin{enumerate}[\rm (a)]
\item $\calV$ is the orthogonal complement of $\calH$.
\item The leaves of the foliation of $\calV$ are totally geodesic.
\end{enumerate}
\end{theorem}
Note that the largest neighborhood $U_x$ in condition \eqref{CondA} is exactly the neighborhood of elements $p \in T^*_xM$ such that both $\exp^{sr}(x,p)$ and $\exp^r(x,\sharp p)$ are well defined. If both $\mathsf{d}^{\tensorh}$ and the distance induced by the Riemannian metric $\tensorg$ are complete metrics, the latter being a sufficient condition for the first, then we may choose $U_x = T^*_xM$, see \cite[Section~7]{Str86}, thus obtaining a global version of Theorem~\ref{th:A}.

\begin{remark}
An equivalent formulation of Theorem~\ref{th:A} is that any curve $\exp^{sr}(x,tp)$ is the horizontal lift of the projection of the curve $\exp^r(x,t\sharp p)$ if and only if (a) and (b) hold. For the definition of horizontal lifts of curves, see Section~\ref{sec:theory}. Another equivalent formulation is that \eqref{CondA} holds if and only if the foliation $\calF = \{ M_b \, \colon \, b \in B\}$ is a totally geodesic Riemannian foliation.
\end{remark}

\subsection{Examples}
\begin{example}[Principal bundles]\label{ex:Montgomery}
Let $G$ be a Lie group with Lie algebra $\mathfrak{g}$, and let $\langle \, \centerdot \, , \,\centerdot \, \rangle_{\mathfrak{g}}$ be a bi-invariant inner product on $\mathfrak{g}$. Let $\pi: M \to B$ be a right principal $G$-bundle over a Riemannian manifold $(B, \widecheck{\tensorg})$. On $M$, for each $A \in \mathfrak{g}$, we have the canonical vector field $\xi_A$ associated to the group action, defined by
$$\xi_A |_x = \left. \frac{d}{dt} (x \cdot \exp^G(tA)) \right|_{t=0}, \qquad x \in M,$$
where $\exp^G$ is the group exponential of $G$. If $\calV = \ker \pi_*$, any element in $\calV_x$ can uniquely be represented as $\xi_A|_x$ for some element $A \in \mathfrak{g}$.

Let $\calH$ be an Ehresmann connection on $\pi$ satisfying $\calH_x \cdot a = \calH_{x \cdot a}$ for any $a \in G$. Then we have a corresponding \emph{connection form} $\omega$, which is a $\mathfrak{g}$-valued one-from uniquely determined by the properties $\ker \omega = \calH$ and $\omega(\xi_A) = A$. Conversely, if $\omega$ is any $\mathfrak{g}$-valued one-from satisfying $\omega(\xi_A) =A$ and $\omega(v \cdot a) = \mathrm{Ad}(a^{-1})\omega(v)$, $\ker \omega$ will always be an Ehresmann connection on $\pi$ invariant under the group action.

Given a connection from $\omega$, we define a Riemannian metric $\tensorg$ on $M$ by
$$\tensorg(v,w) = \widecheck{\tensorg}(\pi_* v, \pi_* w) + \langle \omega(v), \omega(w) \rangle_{\mathfrak{g}}.$$
Let $\nabla^{\tensorg}$ and $\nabla^{\widecheck{\tensorg}}$ denote the Levi-Civita connections for $M$ and $B$, respectively. It can then be verified that for any vector field $\check X$ and $\check Y$ on $B$ and any $A,A_1,A_2 \in \mathfrak{g}$,
\begin{equation} \label{LCpb} 
\begin{cases}
\nabla^{\tensorg}_{h \check X} h \check Y = h \nabla^{\widecheck{\tensorg}}_{\check X} \check Y + \frac{1}{2} \calR(h\check X, h \check Y), \qquad \nabla^{\tensorg}_{\xi_{A_1}} \xi_{A_2} = \frac{1}{2} \xi_{[A_1,A_2]}, \\ \vspace{-0.3cm} \\
\nabla^{\tensorg}_{h \check X} \xi_A = -\nabla^{\tensorg}_{\xi_A} h \check X = -\frac{1}{2} \sharp \tensorg(\xi_A, \calR(h\check X, \, \centerdot \,)),
\end{cases}
\end{equation}
and these relations uniquely determine $\nabla^{\tensorg}$. From \eqref{LCpb} it follows that $\{M_{b}\, \colon \, b \in B\}$ is a totally geodesic Riemannian foliation.

We also have the following two observations.
\begin{enumerate}[\rm (i)]
\item For any vector $v \in TM$, $(\nabla^{\tensorg}_v \omega)(v) =0$, so for any geodesic $\gamma$ in $M$, $\omega(\dot \gamma) = \omega(\pr_{\calV} \dot \gamma)$ is a constant.
\item Since each $M_b$ is a totally geodesic submanifold of $M$ and its metric comes from a bi-invariant metric, we have that for any $v \in \calV_x$, $\exp^r(x, v) = x \cdot \exp^G(\omega(v)).$
\end{enumerate}
We use (i) and (ii) to write the result of Theorem~\ref{th:TGRF} as
\begin{align} \label{gauge}
\exp^{sr}(x,tp) & = \exp^r\left(\exp^r(x,t\sharp p),- t \pr_{\calV} P_t \sharp p\right) \\ \nonumber
& = \exp^r(x,t\sharp p) \cdot \exp^G(- t \omega( P_t \sharp p)) \\ \nonumber
& =  \exp^r(x,t\sharp p) \cdot \exp^G(- t \omega( \sharp p)).
\end{align}
The latter relation was fist observed in \cite[Theorem 11.8]{Mon02}. Projections of the curves in \eqref{gauge} are the trajectories of particles in $B$ with gauge in $\mathfrak{g}^*$ and with Yang-Mills field given by $- \omega(\calR(\, \centerdot \, , \, \centerdot \,))$. For more information, see also \cite{Mon84} or \cite[Chapter 12]{Mon02}. See also \cite[Section~2]{Gro12} for a generalization of this idea to general submersions.

\end{example}

\begin{example}[Octonionic Hopf fibration]
Let us consider the case of a Riemannian submerison $\pi\colon S^m \to B$ with connected totally geodesic fibers, where $S^m$ is the unit sphere with its usual round metric. According to~\cite[Theorem 3.5]{Esc75}, such a Riemannian submersion is necessarily a Hopf fibration whenever $1\leq \dim B \leq m-1$. To be more specific, let ${\mathbb H}$ denote the division algebra of quaternions, and let ${\mathbb C}P^n$ and ${\mathbb H}P^n$ denote the complex and quaternionic projective $n$-spaces, respectively. Let $S^m(r)$ denote the $m$-dimensional sphere of radius $r$. Then any such submersion $\pi:S^m \to B$ is contained in the list
\[
S^1\to S^{2n+1}\to{\mathbb C}P^n,\quad S^3\to S^{4n+3}\to{\mathbb H}P^n,\quad n\geq2,
\]
\[
S^1\to S^3\to S^2(\tfrac12),\quad S^3\to S^7\to S^4(\tfrac12),\quad S^7\to S^{15}\to S^8(\tfrac12).
\]
With the exception of $S^7\to S^{15}\to S^8(\tfrac12)$, all of these submersions can be given structures of principal $\mathrm{U}(1)$- or $\mathrm{SU}(2)$-bundles. The formulas of Example~\ref{ex:Montgomery} have been successfully applied to study the normal sub-Riemannian geodesics for these fibrations listed above, see~\cite{GM12}.

The fibration $S^7\to S^{15}\to S^8(\tfrac12)$ is called the \emph{octonionic Hopf fibration} and requires some more explanation. Here we follow~\cite[Section 6]{GWZ86}. Let ${\mathbb O}$ denote the algebra of octonions with ${\mathbb O}P^1$ being the octonionic projective line. Consider the subsets of ${\mathbb O}\times{\mathbb O}$ given by
\[
L_m=\{(u,mu)\colon u\in{\mathbb O}\}\quad\mbox{for }m\in{\mathbb O},\quad L_\infty=\{(0,u)\colon u\in {\mathbb O}\}.
\]
Given any point $(x,y)\in S^{15}\subset{\mathbb O}\times{\mathbb O}={\mathbb R}^{16}$, we map it to the unique $m\in{\mathbb O}P^1$ such that $(x,y)\in L_m$. The octonionic Hopf fibration has no principal bundle structure, in fact, there are no non-vanishing vertical vector fields, see~\cite[Theorem A]{OPPV13}. However our results, Theorem~\ref{th:TGRF} and Theorem~\ref{th:A}, do hold in this case.


\end{example}

\section{Proofs}\label{sec:proofs}
\subsection{Connections, horizontal and vertical lifts} \label{sec:theory}
Let $\pi: M \to B$ be a submersion with Ehresmann connection~$\calH$ as defined in Section~\ref{sec:submersion}. Then for a given vector $\check v\in T_xB$, the unique element $v\in \calH_x$ such that $\pi_* v = \check v \in T_{\pi(x)}B$ is called \emph{the horizontal lift} of $\check v$. We will write this element as $v = h_x \check v$. Furthermore, if $\check X$ is a vector field on $B$, then $h\check X$ is the vector field on $M$ with values in $\calH$ given by the formula $h\check X|_x := h_x \check X|_{\pi(x)}.$

For a given absolutely continuous curve $\check{\gamma}:[0,T] \to B$, the horizontal lift $\gamma$ of $\check{\gamma}$ to $x \in M_{\check\gamma(0)}$ is the solution of the initial value problem
$$\dot \gamma(t) = h_{\gamma(t)} \dot{\check{\gamma}}, \qquad \gamma(0) = x.$$
This problem clearly has a unique solution, but the horizontal lift $\gamma(t)$ may in general only exist for sufficiently small values of $t$.

Let $\calV = \ker \pi_*$ be the vertical bundle and let $\pr_{\calH}$ and $\pr_{\calV}$ be the respective projections to $\calH$ and $\calV$ with respect to the decomposition $TM = \calH \oplus \calV$. Then \emph{the curvature} $\calR$ of $\calH$ is the vector valued two-form, given by equation
\begin{equation} \label{CurvCalH}
\calR(v,w) = \pr_{\calV} [ \pr_{\calH} X, \pr_{\calH} Y]|_x, \qquad v,w \in T_x M,
\end{equation}
where $X$ and $Y$ are any vector fields satisfying $X|_x = v$ and $Y|_x = w$. It is simple to verify that formula \eqref{CurvCalH} is independent of the choice of vector fields $X$ and $Y$.

We can use the same terminology in a more general setting and define the curvature $\calR$ of $\calH$ by \eqref{CurvCalH} whenever we have some decomposition of the tangent bundle into a direct sum $TM = \calH \oplus \calV$. However, in this case, $\calH$ will also have \emph{a cocurvature} $\overline{\calR}$ analogously given by
$$\overline{\calR}(v,w) = \pr_{\calH} [ \pr_{\calV} X, \pr_{\calV} Y]|_x, \qquad X|_x = v, Y|_x = w.$$
Clearly $\overline{\calR} = 0$ if and only if $\calV$ is integrable. For more information, see~\cite[Chapter III]{KMS93}.

In what follows, we will also need \emph{vertical lifts}, that exist whenever we have a vector bundle $\Pi:\calE \to M$ over a manifold. For any $e_1,e_2 \in \calE_x$, we define the vertical lift of $e_2$ at $e_1$ as
$$\vl_{e_1} e_2 = \left. \frac{d}{dt} ( e_1 + t e_2) \right|_{t=0} \in T_{e_1} \calE.$$
Note that any element of $\ker \Pi_*$ can be written as a vertical lift. Similarly, we can lift a section $\alpha \in \Gamma(\calE)$ of $\calE$ to a vector field $\vl \alpha$ on $\calE$ by $\vl \alpha |_e = \vl_e \alpha|_{\Pi(e)}$. For us, the particular case of $\Pi^M:T^*M \to M$ will be important and the fact that $\ker \Pi_*^M$ are spanned by vertical lifts of one-forms.

\subsection{Connections and the symplectic form}
Let $\vartheta$ be the Liouville one-form on the cotangent bundle $T^*M$ given by $\vartheta|_p(v) =p (\Pi^M_* v)$, and let $\varsigma = -d\vartheta$ be the canonical symplectic form. Consider any affine connection $\nabla$ on $M$ with torsion tensor $T^\nabla$. There is a unique Ehresmann connection $\calE^\nabla$ on $\Pi^M$ such that a smooth curve $\lambda(t)$ in $T^*M$ is tangent to $\calE^{\nabla}$ if and only if $\lambda(t)$ is parallel along $\gamma(t) = \Pi^M(\lambda(t))$, see \cite[Chapter~III]{KMS93}. Then we can write $T(T^*M) = \calE^{\nabla} \oplus (\ker \Pi^M_*)$, where $\calE^{\nabla}$ is spanned by horizontal lifts of vector fields on $M$, while $\ker \Pi^M_*$ is spanned by vertical lifts of forms on~$M$. As a consequence, we can completely describe $\varsigma$ by its values on such elements.

Let $X$ and $Y$ be vector fields on $M$ with horizontal lifts $hX$ and $hY$ and let $\alpha$ and $\beta$ be forms on $M$ with vertical lifts $\vl \alpha$ and $\vl \beta$. Then it is simple to verify from the definition of $\varsigma$ that
\begin{equation} \label{symplectic} \left\{ \begin{array}{rcl} \varsigma(hX|_p , hY|_p) & =& - p(T^{\nabla}(X,Y)), \\  \varsigma(hX|_p, \vl \alpha|_p) & =& \alpha(X|_{\Pi^M(p)}), \\ \varsigma(\vl \alpha|_p, \vl \beta|_p) & =& 0. \end{array} \right. \end{equation}

\subsection{Proof of Proposition~\ref{prop:NormGeo}}
For any (possibly degenerate) cometric $\tensors^* \in \Gamma(\Sym^2 TM)$, define the Hamiltonian $H^{\tensors^*}(p) := \frac{1}{2} \tensors^*(p,p).$ Define the vector $\shs p$ by $\tensors^*(\alpha,p) = \alpha(\shs p) $.
Let $\nabla$ be any connection on $M$ and write $\vec{H}^{\tensors^*} = A + B$ where $A$ and $B$ have values in $\calE^{\nabla}$ and $\ker \Pi^M$, respectively. Then, for any vector field $X$ and one-form $\alpha$ on $M$, we use \eqref{symplectic} to get
\begin{eqnarray*} 
dH^{\tensors^*}(\vl \alpha)|_p & =& \tensors^*(\alpha,p) = \alpha(\shs p) = \varsigma(\vec{H}^{\tensors^*}, \vl \alpha) = \alpha(\Pi^M_* A)|_p, \\
dH^{\tensors^*}(hX)|_p & = & \frac{1}{2} (\nabla_X \tensors^*)(p,p) = \varsigma( \vec{H}^{\tensors^*}, hX)|_p \\
&= & - p(T^\nabla(\Pi^M_* A, X)) - \varsigma(hX, B).
\end{eqnarray*}
It follows that
\begin{equation} \label{HvecS}\vec{H}^{\tensors^*} = h_p \shs p - \vl_p \left(  p T^{\nabla}(\shs p, \centerdot) + \frac{1}{2} (\nabla_{\centerdot} \tensors^*)(p,p) \right).\end{equation}
In order to obtain the result, put $\tensors^* = \tensorh^*$ and use that for any curve $\lambda(t)$ in $T^*M$ with projection $\gamma(t)$, we have
$$\dot \lambda(t) = h_{\lambda(t)} \dot \gamma(t) + \vl_{\lambda(t)} \nabla_{\dot \gamma} \lambda(t).$$

\subsection{Proof of Lemma~\ref{lemma:commute}}
We begin by noting that the torsion $T^{\rnabla}$ of $\rnabla$ is given by $T^{\rnabla} = -\calR - \overline{\calR}$, where $\calR$ and $\overline{\calR}$ denote the curvature and the cocurvature of $\calH$ respectively. From equation \eqref{HvecS}, we get that the Hamiltonian vector fields are given by
\begin{eqnarray}
\vec{H}^{\tensorh^*}|_p &=& h_p \shh p + \vl_p \left(p\calR(\shh p, \, \centerdot \,) + (\rnabla_{\, \centerdot \,} \tensorh^*)(p,p) \right), \nonumber\\ \label{verHam}
\vec{H}^{\tensorv^*}|_p & =&  h_p \shv p + \vl_p \left(p\overline{\calR}(\shv p, \, \centerdot \,) + (\rnabla_{\, \centerdot \,} \tensorv^*)(p,p) \right).
\end{eqnarray}
It now follows that from equation \eqref{symplectic} that
\begin{align*}
   \left\{ H^{\tensorh^*} |_p, H^{\tensorv^*} |_p \right\}  &= \varsigma(\vec{H}^{\tensorh^*}|_p, \vec{H}^{\tensorv^*}|_p) \\ &=  -p T^{\rnabla}(\shh p, \shv p)  +p \overline{\calR}(\shv p, \shh p) + (\rnabla_{\shh p} \tensorv^*)(p,p) \\
&  \quad- p \calR(\shh p, \shv p) - (\rnabla_{\shv p} \tensorh^*)(p,p) \\
& = (\rnabla_{\shh p} \tensorv^*)(p,p) - (\rnabla_{\shv p} \tensorh^*)(p,p)
\end{align*}
Since $\calH$ and $\calV$ are orthogonal with respect to $\tensorg$, we obtain that $\shh p = \pr_{\calH} \sharp p = \sharp \pr_{\calH}^* p$, and similar relations hold for $\shv p$. Then
\begin{align*}
  (\rnabla_{\shh \alpha} \tensorv^*)(\alpha,\alpha)  & = (\pr_{\calH} \sharp \alpha) \tensorg(\sharp \alpha, \pr_{\calV} \sharp \alpha) - 2 (\rnabla_{\pr_{\calH} \sharp \alpha} \alpha) (\pr_{\calV} \sharp \alpha) \\
  & = - (\rnabla_{\pr_{\calH} \sharp \alpha} \tensorg)(\pr_{\calV} \sharp \alpha, \pr_{\calV} \sharp \alpha),
\end{align*}
and similarly $(\rnabla_{\shv \alpha} \tensorh^*)(\alpha, \alpha) = - (\rnabla_{\pr_{\calV} \sharp \alpha} \tensorg)(\pr_{\calH} \sharp \alpha, \pr_{\calH} \sharp \alpha)$. It follows that $\{\vec{H}^{\tensorh} , \vec{H}^{\tensorv} \} = 0$ if and only if the map
\begin{equation} \label{Map} v \mapsto (\rnabla_{\pr_{\calV} v}\tensorg)(\pr_{\calH} v, \pr_{\calH} v) - (\rnabla_{\pr_{\calH} v}\tensorg)(\pr_{\calV} v, \pr_{\calV} v), \quad v \in TM,\end{equation}
vanishes. Notice that the first term in the above map is bilinear in $\pr_{\calH} v$ and linear in $\pr_{\calV} v$ and vice versa for the second term. The map \eqref{Map} is hence zero if and only if both 
\[
(\rnabla_{\pr_{\calH} v} \tensorg)(\pr_{\calV} w, \pr_{\calV} w) = 0 \mbox{ and }(\rnabla_{\pr_{\calV} v} \tensorg)(\pr_{\calH} w, \pr_{\calH} w) = 0 
\]
holds for any $v,w \in TM$. The result now follow from the identity \eqref{CovDerg}.

\subsection{Proof of Theorem~\ref{th:TGRF}} \label{sec:ProofTGRF}
Write $H^{\tensorg} = H^{\tensorh} + H^{\tensorv}$. From Lemma~\ref{lemma:commute}, we know that
$$e^{t\vec{H}^{\tensorh}}(p) = e^{-t\vec{H}^{\tensorv}} \circ e^{t\vec{H}^{\tensorg}}(p),$$
whenever both sides are defined. Now, since $H^{\tensorg}$ is the Hamiltonian of a Riemannian metric, we have that $\alpha(t) = e^{t\vec{H}^{\tensorg}}(p)$ satisfies $\Pi^M(\alpha(t)) = \exp^r(x, t\sharp p)$. Furthermore, if $P_t$ is the parallel transport of vectors with respect to the Levi-Civita connection, then $\sharp \alpha(t) = P_t \sharp p.$

Next, from equation \eqref{verHam}, we know that $\lambda(t) = e^{t\vec{H}^{\tensorv}}(p)$ is a solution to
\begin{equation}\label{eq:fibergeodesic}
\Pi^M(\lambda(t)) = \gamma(t), \qquad \dot \gamma(t) = \shv \lambda(t), \qquad \rnabla_{\dot \gamma} \lambda(t) =0, \quad \lambda(0) = p.
\end{equation}
However, since $\rnabla$ preserves $\Ann(\calH)$ and $\Ann(\calV)$, we might as well consider equation~\eqref{eq:fibergeodesic} with $\lambda(t)$ replaced by $\lambda^{\calV}(t) =\pr_{\calV}^* \lambda(t)$, which satisfies $ \sharp \lambda^{\calV}(0) = \pr_{\calV} \sharp p.$ Finally, since $\lambda(t)$ is a curve in $\Ann(\calH)$, $\gamma(t)$ is tangent to $\calV$ and $\calF$ is a totally geodesic foliation, we have $\rnabla_{\dot \gamma} \lambda(t) = \nabla_{\dot \gamma}^{\tensorg} \lambda(t) = 0$. This completes the proof.

\subsection{Proof of Theorem~\ref{th:A}}
By Proposition~\ref{prop:NormGeo}, if $\gamma(t) = \exp^{sr}(x, t p)$ and $\eta(t) = \exp^r( x, t \sharp p)$, then $\dot \gamma(0) = \shh p$, while $\dot \eta(0) = \sharp p$.

First assume that \eqref{CondA} holds. Let us pick an element $p \in U_x$, such that $\shh p = 0$. Then $\gamma(t) = \exp^{sr}(x, tp)$ is a constant curve, so $\eta(t) = \exp(x,t\sharp p)$ must be a contained in $M_b$ with $b = \pi(x)$. Hence, we must have $\pi_* \dot \eta(0) = \pi_* \sharp p = 0$, so it follows that $\sharp$ maps $\Ann(\calH)$ into $\calV$, and so $\calV$ is orthogonal to~$\calH$. Furthermore, since any geodesic in $M$ which starts in and is tangent to $M_b$ remains in $M_b$, we know that $M_b$ is a totally geodesic submanifold.

The converse statement follows from Theorem~\ref{th:TGRF} and the fact that (a) and (b) imply that the foliation $\{ M_b \, \colon \, b \in B\}$ is a totally geodesic Riemannian foliation.

\end{document}